\documentstyle[amscd,amssymb,verbatim,epsf]{amsart}

\begin{document}
\theoremstyle{plain}
\newtheorem{Thm}{Theorem}
\newtheorem{Cor}{Corollary}
\newtheorem{Con}{Conjecture}
\newtheorem{Main}{Main Theorem}
\newtheorem{Lem}{Lemma}
\newtheorem{Prop}{Proposition}

\theoremstyle{definition}
\newtheorem{Def}{Definition}
\newtheorem{Note}{Note}
\newtheorem{Ex}{Example}

\theoremstyle{remark}
\newtheorem{notation}{Notation}
\renewcommand{\thenotation}{}

\errorcontextlines=0
\numberwithin{equation}{section}
\renewcommand{\rm}{\normalshape}%

\title
   {On Hamilton's Characteristic Functions for Reflection}

\author{Brendan Guilfoyle}
\address{Brendan Guilfoyle\\
          Department of Mathematics and Computing \\
          Institute of Technology, Tralee \\
          Clash \\
          Tralee  \\
          Co. Kerry \\
          Ireland.}
\email{brendan.guilfoyle@@ittralee.ie}

\author{Wilhelm Klingenberg}
\address{Wilhelm Klingenberg\\
 Department of Mathematical Sciences\\
 University of Durham\\
 Durham DH1 3LE\\
 United Kingdom.}
\email{wilhelm.klingenberg@@durham.ac.uk }

\date{January 27, 2006}

\begin{abstract}
We review the complex differential geometry of the space of oriented affine lines in 
${\Bbb{R}}^3$ and give a description of Hamilton's characteristic functions for reflection in an 
oriented C$^1$ surface in terms of this geometry.

\end{abstract}

\maketitle

\section{Introduction}

In a series of classic papers \cite{Ham1} \cite{Ham2} \cite{Ham3}, Hamilton presented his system of geometric optics
whereby the evolution of systems of rays undergoing reflection and refraction at interfaces of homogenous media 
was described, both in general terms and in specific cases. In particular, Hamilton introduced three characteristic 
functions for an optical system by which he could give algebraic expressions for reflection and refraction of a given 
system of rays off simple geometric surfaces, such as spheres or other surfaces of revolution. This he carried out in 
closed form for some systems and in approximation for others. For a beautiful exposition of this work, see Synge's 
book \cite{synge}.

The purpose of this paper is to review the complex differential geometry of the space of oriented affine lines in 
${\Bbb{R}}^3$ and to give a description of Hamilton's characteristic functions for reflection in an 
oriented C$^1$ surface in terms of this geometry. This formalism, developed in \cite{gak2}, has proven 
useful in the scattering of waves off surfaces \cite{gak3}, as well as the computation of caustics \cite{gak4} and 
approximations of the Casimir force \cite{gaks}. 

In the following section we explain the background of this formalism. Section 3 introduces and 
describes the three characteristic functions for reflection in an arbitrary surface in terms of the complex geometry on 
the space of lines (Theorems \ref{t:T} to \ref{t:V}).

\section{The Complex Geometry of Reflection}

\subsection{The Space of Oriented Lines}

The fundamental building block of geometric optics in a homogeneous isotropic medium is a {\it ray}, 
or oriented straight line, along which light 
is postulated to propagate. The set ${\Bbb{L}}$ of all such oriented lines in Euclidean ${\Bbb{R}}^3$ is a 4-dimensional 
manifold that can be identified with the total space of the tangent bundle to the 2-sphere \cite{hitch}. To 
see this, recall that an oriented line is uniquely determined by a pair of vectors: its direction vector $\vec{V}$ and 
the vector $\vec{U}$ from the origin to the point on the line closest to the origin. The former can be taken to be a unit 
vector and is perpendicular to the latter. Thus
\[
{\Bbb{L}}\cong\{\;\;(\vec{U},\vec{V})\in{\Bbb{R}}^3\times{\Bbb{R}}^3\;\;|\quad|\vec{V}|=1\quad\vec{U}\cdot\vec{V}=0\;\;\},
\]
which is clearly homeomorphic to the tangent bundle to the 2-sphere.

Our approach to geometric optics is to study optical systems as submanifolds of the space ${\Bbb{L}}$. In order to 
describe this algebraically we now introduce a local coordinate atlas on  ${\Bbb{L}}$.

Let $\xi$ be the local complex coordinate on the unit 2-sphere in ${\Bbb{R}}^3$
obtained by stereographic projection from the south pole. In terms of the standard spherical polar angles $(\theta,\phi)$
about the north pole, we have $\xi=\tan(\frac{\theta}{2})e^{i\phi}$. 

This can be extended to complex
coordinates $(\xi,\eta)$ on ${\Bbb{L}}$ minus the tangent space over the south
pole, as follows. Note that a tangent vector $\vec{\mbox{X}}$ to the 2-sphere can  
always be expressed as a linear combination:
\[
\vec{\mbox{X}}=\eta\frac{\partial}{\partial\xi}+\bar{\eta}\frac{\partial}{\partial\bar{\xi}},
\]
for a complex number $\eta$. We can thus identify the real tangent vector 
$\vec{\mbox{X}}$ on the 2-sphere 
(and hence the ray in ${\Bbb{R}}^3$) with the two complex numbers ($\xi,\eta$). Loosely speaking,
$\xi$ determines the direction of the ray, and $\eta$ determines its perpendicular distance
vector to the origin.

The coordinates ($\xi,\eta$) do not cover all of ${\Bbb{L}}$ - they omit all of the lines 
pointing directly downwards. However, the construction can also be carried out 
using stereographic projection from the north pole, yielding a coordinate system that 
covers all of ${\Bbb{L}}$ except for the lines pointing directly upwards. The whole of the space of oriented lines is 
covered by these two coordinate patches  and the transition functions between
them are holomorphic. In what follows we work in the patch 
that omits the south direction.

\begin{Def}
The canonical bundle map $\pi:{\Bbb{L}}\rightarrow{\Bbb{P}}^1$ is given by mapping an oriented line to its
direction. Here, and throughout, we write ${\Bbb{P}}^1$ for S$^2$ with the standard complex structure, and
denote the composition of projections ${\Bbb{L}}\times {\Bbb{R}}\rightarrow{\Bbb{L}}\rightarrow{\Bbb{P}}^1$ by $\pi_1$.
\end{Def}

One of the key relationships between ${\Bbb{L}}$ and ${\Bbb{R}}^3$ is the map 
$\Phi:{\Bbb{L}}\times{\Bbb{R}}\rightarrow{\Bbb{R}}^3$. 

\begin{Def}
The map $\Phi$ takes 
$((\xi,\eta),r)\in{\Bbb{L}}\times{\Bbb{R}}$ to the point on the oriented line ($\xi,\eta$) in ${\Bbb{R}}^3$ 
that lies an affine parameter distance $r$
from the point on ($\xi,\eta$) closest to the origin. 
\end{Def}

If $\Phi((\xi,\eta),r)=(z(\xi,\eta,r),t(\xi,\eta,r))$, where $z=x^1+ix^2$, $t=x^3$ and ($x^1$, $x^2$, $x^3$) are Euclidean
coordinates in ${\Bbb{R}}^3$, then we have the 
following coordinate expressions \cite{gak2}:
\begin{equation}\label{e:coord}
z=\frac{2(\eta-\overline{\eta}\xi^2)+2\xi(1+\xi\overline{\xi})r}{(1+\xi\overline{\xi})^2}
\qquad\qquad
t=\frac{-2(\eta\overline{\xi}+\overline{\eta}\xi)+(1-\xi^2\overline{\xi}^2)r}{(1+\xi\overline{\xi})^2}.
\end{equation}
This map is of crucial importance when describing surfaces in ${\Bbb{R}}^3$, 
as we explain below. It is also useful to invert these relationships:
\begin{equation}\label{e:coord2a}
\eta={\textstyle{\frac{1}{2}}}\left(z-2t\xi-\bar{z}\xi^2\right)
\qquad\qquad
r=\frac{\xi\bar{z}+\bar{\xi}z+(1-\xi\bar{\xi})t}{1+\xi\bar{\xi}}.
\end{equation}
These give the perpendicular distance vector from the origin and the distance $r$ for a line through the point 
$(z,t)\in{\Bbb{R}}^3$ with direction $\xi\in{\Bbb{P}}^1$.

For later use we introduce the following definition:

\begin{Def}
For a subset $K\subset{\Bbb{R}}^3$, the subset ${\Bbb{U}}_K\subset{\Bbb{L}}$ is defined to be the 
set of oriented lines that intersect $K$:
\[
{\Bbb{U}}_K\equiv\left\{\;\gamma\in{\Bbb{L}}\;\left|\;\Phi(\gamma,\cdot)\cap K\neq\emptyset\;\right.\right\}.
\]
\end{Def}

\begin{Ex}
For $p\in{\Bbb{R}}^3$ the set ${\Bbb{U}}_p$ is the set of oriented lines through $p$, which is a sphere in ${\Bbb{L}}$. 
In terms of the coordinates above, if $p$ is ($z$,$t$), then
\[
{\Bbb{U}}_{(z,t)}\equiv\{\;(\xi,\eta)\in{\Bbb{L}}\;|\;\eta={\scriptstyle{\frac{1}{2}}}
             \left(z-2t\xi-\bar{z}\xi^2\right)\;\}.
\]
\end{Ex}

\begin{Ex}
For a C$^1$ surface $S\subset{\Bbb{R}}^3$, ${\Bbb{U}}_S$ is a 4-dimensional subset of ${\Bbb{L}}$ whose boundary is
\[
\partial{\Bbb{U}}_S=\{\;\gamma\in{\Bbb{L}}\;|\;\Phi(\gamma,\cdot) {\mbox{ tangent to }}S\;\}.
\] 
\end{Ex}

In its simplest form, an optical system is a 2-parameter family of rays. We refer to a 2-parameter family of
oriented lines in ${\Bbb{R}}^3$ as a {\it line congruence}, which, from our perspective, is a 
surface $\Sigma$ in ${\Bbb{L}}$. For example, a point source corresponds to the 2-parameter family of 
oriented lines that contain the source point, which thus defines a 2-sphere in ${\Bbb{L}}$ - the sphere in Example
1.

For computational purposes, we now give explicit local parameterizations of the line congruence. In practice,
this will be given locally by a map ${\Bbb{C}}\rightarrow{\Bbb{L}}:\mu\mapsto(\xi(\mu,\bar{\mu}),\eta(\mu,\bar{\mu}))$. A
convenient choice of parameterisation will often depend upon the specifics of the situation, but 
our formalism holds for arbitrary parameterisations.

We now describe how to construct surfaces in ${\Bbb{R}}^3$ using line
congruences. Given a line congruence $\Sigma\subset{\Bbb{L}}$, a map
$r:\Sigma\rightarrow {\Bbb{R}}$ determines a map
$\Sigma\rightarrow{\Bbb{R}}^3$ by
$(\xi,\eta)\mapsto\Phi((\xi,\eta),r(\xi,\eta))$ for $(\xi,\eta)\in\Sigma$.
With a local parameterization $\mu$ of $\Sigma$, composition with the above map yields a
map ${\Bbb{C}}\rightarrow{\Bbb{R}}^3$ which comes from 
substituting $\xi=\xi(\mu,\bar{\mu})$, $\eta=\eta(\mu,\bar{\mu})$ and $r=r(\mu,\bar{\mu})$ in equations (\ref{e:coord}). In other words, we choose one point on each line of the 
congruence, thus forming a surface.

\vspace{0.1in}
\setlength{\epsfxsize}{3.5in}
\begin{center}
   \mbox{\epsfbox{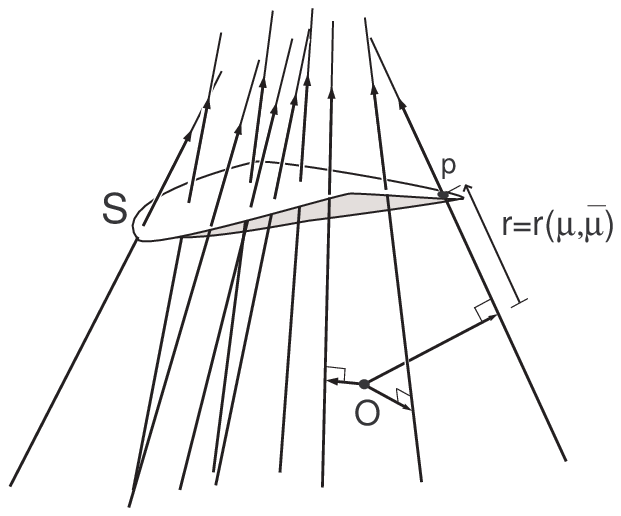}}
\end{center}
\vspace{0.1in}

\begin{Def}
An oriented C$^1$ surface $S$ in ${\Bbb{R}}^3$ gives rise to the map 
$\Phi_S^{-1}:S\rightarrow{\Bbb{L}}\times{\Bbb{R}}$ which takes a point  $p\in S$ to the oriented normal line to $S$ 
crossed with the distance of $p$ to the point on the normal lying closest to the origin. 
\end{Def}

The dual picture of light propagation is to consider the wavefronts, or surfaces that are orthogonal to a
given set of rays. However, not every line congruence has such orthogonal surfaces - indeed, most don't. We refer
to line congruences for which orthogonal surfaces exist as {\it normal congruences}.

\begin{Thm}\cite{gak2}
A parameterised line congruence $(\xi(\mu,\bar{\mu}),\eta(\mu,\bar{\mu}))$ is orthogonal to a surface in ${\Bbb{R}}^3$ iff 
there exists a real function $r(\mu,\bar{\mu})$ satisfying:
\begin{equation}\label{e:intsur}
\partial r=\frac{2\bar{\eta}\partial\xi+2\eta\partial\bar{\xi}}{(1+\xi\bar{\xi})^2},
\end{equation}
where $\partial=\frac{\partial}{\partial\mu}$.

If there exists one solution, there exists a 1-parameter family generated by the real constant of integration.
An explicit parameterization of these surfaces in ${\Bbb{R}}^3$ is given by inserting 
$(\xi(\mu,\bar{\mu}),\eta(\mu,\bar{\mu}))$ and $r=r(\mu,\bar{\mu})$ in  (\ref{e:coord}). 
\end{Thm}

Given a surface $S$ with coordinates $\xi=\xi_0(\mu,\bar{\mu})$, $\eta=\eta_0(\mu,\bar{\mu})$ and $r=r_0(\mu,\bar{\mu})$,
the set ${\Bbb{U}}_S$ is \cite{gak3}:
\begin{equation}\label{e:key1}
{\Bbb{U}}_S=\left\{\;(\xi,\eta)\in{\Bbb{L}}\;\left|\;\eta={\textstyle \frac{(1+\bar{\xi}_0\xi)^2}{(1+\xi_0\bar{\xi}_0)^2}}\eta_0
-{\textstyle \frac{(\xi_0-\xi)^2}{(1+\xi_0\bar{\xi}_0)^2}}\bar{\eta}_0
+{\textstyle \frac{(\xi_0-\xi)(1+\bar{\xi}_0\xi)}{1+\xi_0\bar{\xi}_0}}r_0\;\right.\right\}.
\end{equation}

\subsection{Reflection}

Given an oriented C$^1$ surface $S$ in ${\Bbb{R}}^3$, reflection in $S$ determines a map ${\cal{R}}_S:{\Bbb{U}}_S\rightarrow{\Bbb{U}}_S$.

\vspace{0.1in}
\setlength{\epsfxsize}{3.8in}
\begin{center}
   \mbox{\epsfbox{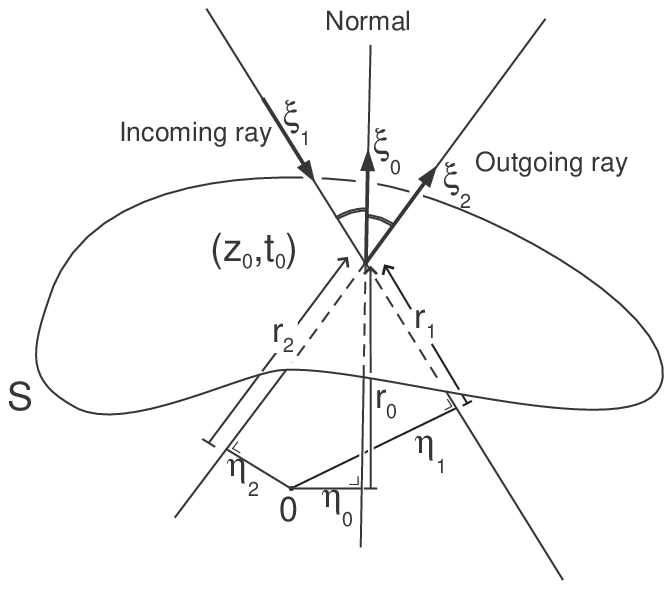}}
\end{center}
\vspace{0.1in}

The algebraic expression for this map is:

\begin{Thm}\label{t:refl} \cite{gak3}
Let $S$ be an oriented C$^1$ surface $S$ in ${\Bbb{R}}^3$ with parameterised normal line
congruence $\xi=\xi_0(\mu,\bar{\mu})$, $\eta=\eta_0(\mu,\bar{\mu})$ and $r=r_0(\mu,\bar{\mu})$ satisfying 
(\ref{e:intsur}). Then  ${\cal{R}}_S:{\Bbb{U}}_S\rightarrow{\Bbb{U}}_S$ takes ($\xi_1,\eta_1$) to 
${\cal{R}}_S(\xi_1,\eta_1)=(\xi_2,\eta_2)$: 
\begin{equation}\label{e:reflawa}
\xi_2=\frac{2\xi_0\bar{\xi}_1+1-\xi_0\bar{\xi}_0}
           {(1-\xi_0\bar{\xi}_0)\bar{\xi}_1-2\bar{\xi}_0},
\end{equation}
\begin{equation}\label{e:key2}
\eta_2={\textstyle \frac{(\bar{\xi}_0-\bar{\xi}_1)^2}
         {((1-\xi_0\bar{\xi}_0)\bar{\xi}_1-2\bar{\xi}_0)^2}}\eta_0
       -{\textstyle\frac{(1+\xi_0\bar{\xi}_1)^2}
         {((1-\xi_0\bar{\xi}_0)\bar{\xi}_1-2\bar{\xi}_0)^2}}\bar{\eta}_0
+{\textstyle\frac{(\bar{\xi}_0-\bar{\xi}_1)(1+\xi_0\bar{\xi}_1)(1+\xi_0\bar{\xi}_0)}
         {((1-\xi_0\bar{\xi}_0)\bar{\xi}_1-2\bar{\xi}_0)^2}}r_0,
\end{equation}
where $(\xi_0,\eta_0,r_0)\in\Phi_S^{-1}(S)$.

If $r_1$ and $r_2$ are the distances of the point of reflection from
the point closest to the origin on the incoming and reflected rays,
respectively, then
\begin{equation}\label{e:rs}
r_2=r_1+\frac{2(|\xi_0-\xi_1|^2-|1+\bar{\xi}_0\xi_1|^2)}
    {(1+\xi_0\bar{\xi}_0)(1+\xi_1\bar{\xi}_1)}r_0.
\end{equation}

\end{Thm}

\section{Hamilton's Characteristic Functions For Reflection in a Surface}

Consider a ray ($\xi_1,\eta_1$) entering an optical instrument and emerging, after a number of 
reflections and refractions, as the ray ($\xi_2,\eta_2$). Let $p_1$ and $p_2$ be arbitrary points on the incoming
and outgoing rays, $q_1$ and $q_2$ be the points on the incoming and outgoing rays that lie closest to the origin 
(respectively). 

\vspace{0.1in}
\setlength{\epsfxsize}{5.0in}
\begin{center}
   \mbox{\epsfbox{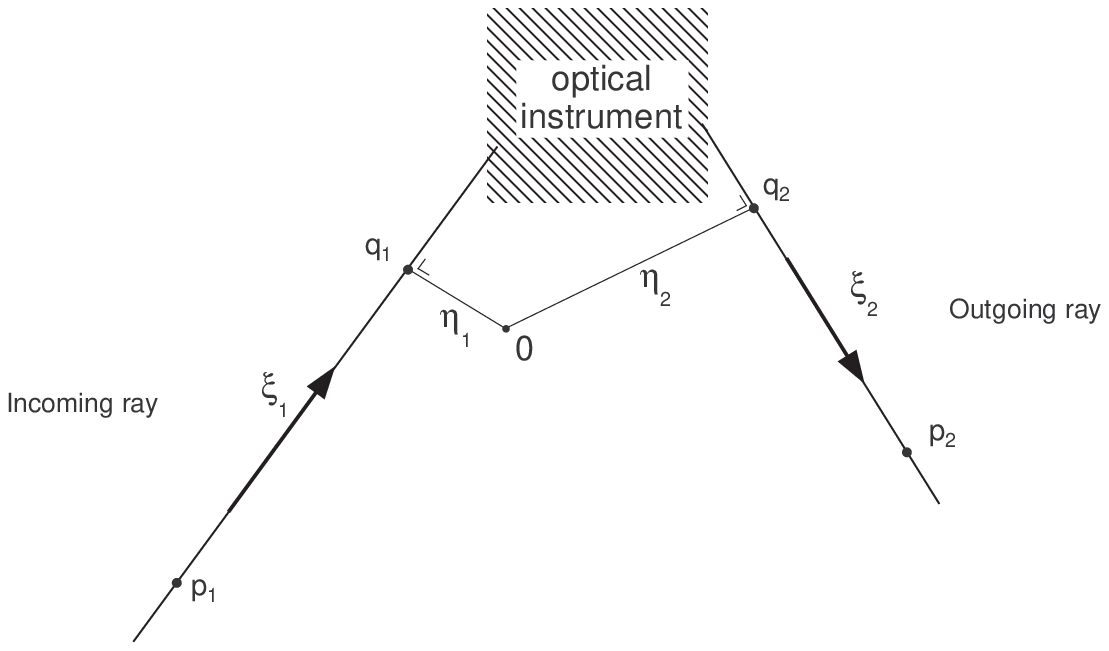}}
\end{center}
\vspace{0.1in}

Hamilton defined three different functions, which he termed {\it characteristic functions} and denoted $T$, $V$ and $W$.

\begin{Def}
The function $T$, which is referred to as the {\it angle characteristic function}, takes the directions $\xi_1$ and $\xi_2$
to the distance along the rays between $q_1$ and $q_2$. 
\end{Def}
\begin{Def}
The function $W$, which is referred to as the {\it mixed 
characteristic function}, takes the point $p_1$ and the direction $\xi_2$ to the distance along the rays between 
$p_1$ and $q_2$. 
\end{Def}
\begin{Def}
The function $V$, which is referred to as the {\it point characteristic function}, takes 
the points $p_1$ and $p_2$ to the distance along the rays between $p_1$ and $p_2$. 
\end{Def}

The function $V$ is, of course, commonly called the Hamiltonian of the optical system.
 
In the above, while the points $q_i$ are on the initial or final rays,
they may be on the opposite side of the initial or final interface and
so lie at points
which the light may never physically pass through. However, as long as
the associated distances are counted negatively, the functions are
well defined and have all of the formal properties required for our purposes.

The exact domain of definition of these functions will depend upon the optical instrument under consideration. In 
addition, the characteristic functions may be multivalued, since there may be more than one sequence of rays with the
properties required. 

The domain of $T$ is obviously a subset of ${\Bbb{P}}^1\times{\Bbb{P}}^1$, while the domain of $V$ 
is a subset of ${\Bbb{R}}^3\times{\Bbb{P}}^1$ and the domain of $W$ is a subset of ${\Bbb{R}}^3\times{\Bbb{R}}^3$. More
specifically, for reflection in a surface:

\begin{Prop}
Let $S$ be an oriented C$^1$ immersed surface in ${\Bbb{R}}^3$, then for reflection in $S$:
\[
{\mbox{Dom T}}=\left\{\; (\xi_1,\xi_2)\in{\Bbb{P}}^1\times{\Bbb{P}}^1\;\left|\; 
                    {\cal{R}}_S\left(\pi^{-1}(\xi_1)\cap{\Bbb{U}}_S\right)\cap\pi^{-1}(\xi_2)\neq\emptyset\;\right.\right\}
\]
\[
{\mbox{Dom W}}=\left\{\; (p_1,\xi_2)\in{\Bbb{R}}^3\times{\Bbb{P}}^1\;\left|\;  
                   {\cal{R}}_S\left({\Bbb{U}}_{p_1}\cap{\Bbb{U}}_S\right)\cap\pi^{-1}(\xi_2)\neq\emptyset\;\right.\right\}
\]
\[
{\mbox{Dom V}}=\left\{\; (p_1,p_2)\in{\Bbb{R}}^3\times{\Bbb{R}}^3\;\left|\; 
                   {\cal{R}}_S\left({\Bbb{U}}_{p_1}\cap{\Bbb{U}}_S\right)\cap{\Bbb{U}}_{p_2}\neq\emptyset\;\right.\right\}.
\]

\end{Prop}
\begin{pf}
A pair of directions $\xi_1,\xi_2\in{\Bbb{P}}^1$ lie in the domain of $T$ iff there exists a line with direction 
$\xi_1$ which, when reflected in $S$, has direction $\xi_2$. Now $\pi^{-1}(\xi_1)$ is the set of oriented lines with 
direction $\xi_1$, while $\pi^{-1}(\xi_1)\cap{\Bbb{U}}_S$ is the set of oriented lines with direction $\xi_1$ that 
intersect $S$. Thus, reflecting this set gives ${\cal{R}}_S\left(\pi^{-1}(\xi_1)\cap{\Bbb{U}}_S\right)$ and, for 
$(\xi_1,\xi_2)$ to be in the domain of $T$, the reflected set must contain a line with direction $\xi_2$. In other words,
${\cal{R}}_S\left(\pi^{-1}(\xi_1)\cap{\Bbb{U}}_S\right)\cap\pi^{-1}(\xi_2)\neq\emptyset$, as claimed.

For $W$, $(p_1,\xi_2)\in{\Bbb{R}}^3\times{\Bbb{P}}^1$ is in the domain of $W$ iff there exists a line through $p_1$ whose 
reflection in $S$ has direction $\xi_2$. The oriented lines through $p_1$ that intersect $S$ are 
${\Bbb{U}}_{p_1}\cap{\Bbb{U}}_S$. Thus, for $(p_1,\xi_2)$ to be in the domain of $W$, the reflection of this set, 
${\cal{R}}_S\left({\Bbb{U}}_{p_1}\cap{\Bbb{U}}_S\right)$ must contain
a line with direction $\xi_2$. That is, we must have 
${\cal{R}}_S\left({\Bbb{U}}_{p_1}\cap{\Bbb{U}}_S\right)\cap\pi^{-1}(\xi_2)\neq\emptyset$, as claimed.

Finally, the set of oriented lines through $p_1$ that intersect $S$ is ${\Bbb{U}}_{p_1}\cap{\Bbb{U}}_S$. Thus, for 
$(p_1,p_2)$ to be in the domain of $V$, the reflection
of this set must contain an oriented line that passes through $p_2$. That is,
${\cal{R}}_S\left({\Bbb{U}}_{p_1}\cap{\Bbb{U}}_S\right)\cap{\Bbb{U}}_{p_2}\neq\emptyset$.
\end{pf}

The following propositions describe the domains of the characteristic functions in terms of our coordinates. We assume
that the oriented C$^1$ surface $S$ has normal congruence parameterised by 
$\mu\mapsto(\xi_0(\mu,\bar{\mu}),\eta_0(\mu,\bar{\mu}))$ and $r=r_0(\mu,\bar{\mu})$.

\begin{Prop}
Given $\xi_1,\xi_2\in{\Bbb{P}}^1$, define $\xi_0\in{\Bbb{P}}^1$ by
\begin{equation}\label{e:xi0}
\xi_0=\frac{\xi_1\bar{\xi}_1-\xi_2\bar{\xi}_2+|\xi_1-\xi_2|\;
         [(1+\xi_1\bar{\xi}_1)(1+\xi_2\bar{\xi}_2)]^{\scriptstyle{\frac{1}{2}}}}
         {\bar{\xi}_1(1+\xi_2\bar{\xi}_2)-\bar{\xi}_2(1+\xi_1\bar{\xi}_1)}.
\end{equation}
Then the domain of $T$ is given by:
\begin{equation}\label{e:T1}
{\mbox{Dom T}}=\left\{\; (\xi_1,\xi_2)\in{\Bbb{P}}^1\times{\Bbb{P}}^1\;\left|\; 
                    \xi_0\in\pi_1\circ\Phi_S^{-1}(S)\;\right.\right\}.
\end{equation}
\end{Prop}
\begin{pf}
Given two directions $\xi_1,\xi_2\in{\Bbb{P}}^1$, the direction $\xi_0$ through which $\xi_1$ must be
reflected to obtain $\xi_2$ (or vice versa) is found by solving equation (\ref{e:reflawa}) for $\xi_0$.
The result is (\ref{e:xi0}), and the question then reduces to whether the normal to $S$ has direction $\xi_0$ at
any point. Now, the map $\pi_1\circ\Phi_S^{-1}:S\rightarrow{\Bbb{P}}^1$ takes a point on an oriented C$^1$ surface $S$ to 
the direction of the normal at the point (i.e. the classical Gauss map of $S$). Thus, for $(\xi_1,\xi_2)$ to be
in the domain of $T$ we must have $\xi_0\in\pi_1\circ\Phi_S^{-1}(S)$.
\end{pf}

\begin{Prop}\label{p:W}
The domain of $W$ is obtained as follows. Given $p_1=(z_1,t_1)\in{\Bbb{R}}^3$ and $\xi_2\in{\Bbb{P}}^1$, solve
(if possible) the single complex equation for $\mu$:
\begin{align}\label{e:w1}
z_1[(1-\xi_0\bar{\xi}_0)\bar{\xi}_2-2\bar{\xi}_0]^2
     &-2t_1[(1-\xi_0\bar{\xi}_0)\bar{\xi}_2-2\bar{\xi}_0][2\xi_0\bar{\xi}_2+1-\xi_0\bar{\xi}_0]\nonumber \\
&\qquad -\bar{z}_1[2\xi_0\bar{\xi}_2+1-\xi_0\bar{\xi}_0]^2\nonumber \\
=2(1+\xi_0\bar{\xi}_0)\left((\bar{\xi}_0-\bar{\xi}_2)^2\eta_0\right.&\left.-(1+\xi_0\bar{\xi}_2)^2\bar{\eta}_0
 +(\bar{\xi}_0-\bar{\xi}_2)(1+\xi_0\bar{\xi}_2)(1+\xi_0\bar{\xi}_0)r_0\right).
\end{align}
Then $(p_1,\xi_2)$ is in the domain of $W$ iff $\xi_0(\mu,\bar{\mu})\in\pi_1\circ\Phi_S^{-1}(S)$.
\end{Prop}
\begin{pf}
An oriented line $(\xi_1,\eta_1)$ passes through the point $(z_1,t_1)$ iff ({\it cf.} equation (\ref{e:coord2a}))
\[
\eta_1={\textstyle{\frac{1}{2}}}\left(z_1-2t_1\xi_1-\bar{z}_1\xi_1^2\right),
\]
while inverting the reflection law (\ref{e:reflawa}) gives
\[
\xi_1=\frac{2\xi_0\bar{\xi}_2+1-\xi_0\bar{\xi}_0}
           {(1-\xi_0\bar{\xi}_0)\bar{\xi}_2-2\bar{\xi}_0}.
\]
For the incoming ray to intersect the surface $S$ we must have ({\it cf.} (\ref{e:key1}))
\[
\eta_1={\textstyle \frac{(1+\bar{\xi}_0\xi_1)^2}{(1+\xi_0\bar{\xi}_0)^2}}\eta_0
-{\textstyle \frac{(\xi_0-\xi_1)^2}{(1+\xi_0\bar{\xi}_0)^2}}\bar{\eta}_0
+{\textstyle \frac{(\xi_0-\xi_1)(1+\bar{\xi}_0\xi_1)}{1+\xi_0\bar{\xi}_0}}r_0.
\]
Thus, for an oriented line through $p_1$ to have reflected direction $\xi_2$, the equation to be solved is
obtained by substituting the first two equations above in the last one. The result is as claimed.
\end{pf}

\begin{Prop}\label{p:V}
The domain of $V$ is determined as follows. Given two points $p_1=(z_1,t_1)$ and $p_2=(z_2,t_2)\in{\Bbb{R}}^3$, solve
(if possible) the following two complex equations for $\mu$ and $\xi_1$:
\begin{align}\label{e:domv1}
z_2[(1-\xi_0\bar{\xi}_0)\bar{\xi}_1-2\bar{\xi}_0]^2
     &-2t_2[(1-\xi_0\bar{\xi}_0)\bar{\xi}_1-2\bar{\xi}_0][2\xi_0\bar{\xi}_1+1-\xi_0\bar{\xi}_0]\nonumber \\
+\bar{z}_2[2\xi_0\bar{\xi}_1+1-\xi_0\bar{\xi}_0]^2&=-(1+\xi_0\bar{\xi}_0)^2(\bar{z}_1-2t_1\bar{\xi}_1-z_1\bar{\xi}_1^2)
      \nonumber\\
&\qquad+4(\bar{\xi}_0-\bar{\xi}_1)(1+\xi_0\bar{\xi}_1)(1+\xi_0\bar{\xi}_0)r_0
\end{align}
and
\begin{align}\label{e:domv2}
(1+\xi_0\bar{\xi}_0)^2\left(z_1-2t_1\xi_1+\bar{z}_1\xi_1^2\right)
&=2(1+\bar{\xi}_0\xi_1)^2\eta_0-2(\xi_0-\xi_1)^2\bar{\eta}_0 \nonumber\\
& -2(1+\bar{\xi}_0\xi_1)(\xi_0-\xi_1)(1+\xi_0\bar{\xi}_0)r_0.
\end{align}
Then $(p_1,p_2)$ is in the domain of $V$ iff $\xi_0(\mu,\bar{\mu})\in\pi_1\circ\Phi_S^{-1}(S)$.
\end{Prop}
\begin{pf}
As in the previous proposition, we must have 
\[
\eta_i={\textstyle{\frac{1}{2}}}\left(z_i-2t_i\xi_i-\bar{z}_i\xi_i^2\right),
\]
this time for $i=1,2$. Substituting these into equations (\ref{e:key1}) and (\ref{e:key2}), with the aid
of (\ref{e:reflawa}), yields the two equations in the proposition.
\end{pf}

We now compute the characteristic functions for reflection in an oriented C$^1$ surface $S$. Assume that 
the normal line congruence to $S$ is parameterized by $\xi=\xi_0(\mu,\bar{\mu})$,
$\eta=\eta_0(\mu,\bar{\mu})$ and $r=r_0(\mu,\bar{\mu})$.

The angle characteristic function for reflection in a surface $S$ is given by:

\begin{Thm}\label{t:T}
Hamilton's angle characteristic function for reflection in $S$ is
\[
T(\xi_1,\xi_2)=\pm\frac{2|\xi_1-\xi_2|}{\left[(1+\xi_1\bar{\xi}_1)
    (1+\xi_2\bar{\xi}_2)\right]^{\scriptstyle{\frac{1}{2}}}}r_0,
\]
where $\xi_0$ is given by (\ref{e:xi0}) and  $(\xi_0,\eta_0,r_0)\in \Phi_S^{-1}(S)$.
\end{Thm}
\begin{pf}
Given $\xi_1$ and $\xi_2$, from the reflection law (\ref{e:reflawa}) we
solve for $\xi_0$ and get the direction $\xi_0$ of the normal to the surface at the point of reflection - equation 
(\ref{e:xi0}). It is 
possible that the surface $S$ never has this normal direction, in which case the function $T$ is not defined at
($\xi_1$,$\xi_2$). On the other hand, there may be more than one point on the surface with normal direction $\xi_0$, in
which case $T$ will be multivalued. Assume that it is defined and that the point of reflection is ($z_0,t_0$) which lies 
at $r=r_0$ along the normal ray ($\xi_0$,$\eta_0$). 

As before, let $r_1$ and $r_2$ be the distances of the point of
reflection from the point closest to the origin on the incoming and
reflected rays, respectively. Then $T(\xi_1,\xi_2)=\pm|r_1-r_2|$ and
by (\ref{e:rs})
\[
T=\pm\frac{2(|1+\bar{\xi}_0\xi_1|^2-|\xi_0-\xi_1|^2)}
    {(1+\xi_0\bar{\xi}_0)(1+\xi_1\bar{\xi}_1)}r_0. 
\]

Substituting equation (\ref{e:xi0}) in this yields the result after some
simplification.
\end{pf}

\vspace{0.2in}

On the other hand the mixed characteristic function for reflection is determined by: 

\begin{Thm}\label{t:W}
Hamilton's mixed characteristic function for reflection in $S$ is:
\[
W((z_1,t_1),\xi_2)=\pm\left|\frac{\bar{\xi}_1z_1+\xi_1\bar{z}_1+(1-\xi_1\bar{\xi}_1)t_1}{1+\xi_1\bar{\xi}_1}
       + \frac{2|\xi_1-\xi_2|}{\left[(1+\xi_1\bar{\xi}_1)
    (1+\xi_2\bar{\xi}_2)\right]^{\scriptstyle{\frac{1}{2}}}}r_0 \right|,
\] 
where 
\begin{equation}\label{e:w4}
\xi_1=\frac{2\xi_0\bar{\xi}_2+1-\xi_0\bar{\xi}_0}
           {(1-\xi_0\bar{\xi}_0)\bar{\xi}_2-2\bar{\xi}_0},
\end{equation}
and $\xi_0\in \pi_1\circ\Phi_S^{-1}(S)$ is a solution of (\ref{e:w1}).
\end{Thm}
\begin{pf}
From Proposition \ref{p:W}, $(z_1,t_1),\xi_2)$ is in the domain of $W$ iff there exists a direction
$\xi_0\in \pi_1\circ\Phi_S^{-1}(S)$ that solves equation (\ref{e:w1}). Assuming that such exists, equation
(\ref{e:w4}) follows from inverting the reflection law (\ref{e:reflawa}).

It is clear that
\[
W((z_1,t_1),\xi_2)=\pm|s_1+r_1-r_2|,
\]
where $s_1$ is the distance from $p_1$ to $q_1$. From the second equation of (\ref{e:coord2a}), we have that
\[
s_1=\frac{\bar{\xi}_1z_1+\xi_1\bar{z}_1+(1-\xi_1\bar{\xi}_1)t_1}{1+\xi_1\bar{\xi}_1},
\]
while from the proof of the previous theorem, we have that
\[
r_1-r_2=\frac{2|\xi_1-\xi_2|}{\left[(1+\xi_1\bar{\xi}_1)
    (1+\xi_2\bar{\xi}_2)\right]^{\scriptstyle{\frac{1}{2}}}}r_0.
\]
The result follows.
\end{pf}

\vspace{0.2in}

Finally, the point characteristic function for reflection can be computed using:

\begin{Thm}\label{t:V}
Hamilton's point characteristic function for reflection in $S$ is:
\[
V((z_1,t_1),(z_2,t_2))=\pm\left|s_1-s_2+ \frac{2|\xi_1-\xi_2|}{\left[(1+\xi_1\bar{\xi}_1)
    (1+\xi_2\bar{\xi}_2)\right]^{\scriptstyle{\frac{1}{2}}}}r_0\right|,
\] 
where  
\begin{equation}\label{e:v4}
s_i=\frac{\bar{\xi}_iz_i+\xi_i\bar{z}_i+(1-\xi_i\bar{\xi}_i)t_i}{1+\xi_i\bar{\xi}_i},
\end{equation}
for $i=1,2$,
\begin{equation}\label{e:e3}
\xi_2=\frac{2\xi_0\bar{\xi}_1+1-\xi_0\bar{\xi}_0}
           {(1-\xi_0\bar{\xi}_0)\bar{\xi}_1-2\bar{\xi}_0},
\end{equation}
and $\xi_0$ and $\xi_1$ solve equations (\ref{e:domv1}) and (\ref{e:domv2}), 
and $\xi_0\in\pi_1\circ\Phi_S^{-1}(S)$.
\end{Thm}
\begin{pf}

From Proposition \ref{p:V}, $(z_1,t_1),(z_2,t_2))$ is in the domain of $V$ iff there exist directions $\xi_0$ and $\xi_1$
that solve equations (\ref{e:domv1}) and (\ref{e:domv2}), and $\xi_0\in\pi_1\circ\Phi_S^{-1}(S)$. 
Assuming that such exists, equation (\ref{e:e3} is just the reflection law (\ref{e:reflawa}), and 
\[
V((z_1,t_1),(z_2,t_2))=\pm\left|s_1-s_2+r_1-r_2\right|.
\]
Finally, the distance from $p_i$ to $q_i$ is given by (\ref{e:v4}).
\end{pf}

\end{document}